\pdfoutput=1
\documentclass{amsart}
\usepackage{graphicx}
\usepackage{amssymb}

\newcommand{\ri}{\operatorname{RI}}
\newcommand{\rii}{\operatorname{RI\!I}}
\newcommand{\riii}{\operatorname{RI\!I\!I}}

\theoremstyle{plain}
\newtheorem{theorem}{Theorem}

\newtheorem{ostlund}{\"{O}stlund Question}

\newtheorem{haggeyazinski}{Hagge-Yazinski Theorem}

\theoremstyle{remark}
\newtheorem{remark}{Remark}

\theoremstyle{definition}
\newtheorem{definition}{Definition}

\begin{document}

\title{On a nontrivial knot projection under (1, 3) homotopy}
\author{Noboru Ito}
\address{Waseda Institute for Advanced Study, 1-6-1 Nishi Waseda, Shinjuku-ku, Tokyo, 169-8050, Japan}
\email{noboru@moegi.waseda.jp}
\address{(Current address) Graduate School of Mathematical Sciences, The University of Tokyo, 3-8-1, Komaba, Meguro-ku, Tokyo, 153-8914, Japan}
\email[(Current)]{noboru@ms.u-tokyo.ac.jp}
\author{Yusuke Takimura}
\address{Gakushuin Boys' Junior High School, 1-5-1, Mejiro, Toshima-ku, Tokyo, 171-0031, Japan}
\email{Yusuke.Takimura@gakushuin.ac.jp}
%\date{\today}
\maketitle

\begin{abstract}
In 2001, \"{O}stlund formulated the question: are Reidemeister moves of types 1 and 3 sufficient to describe a homotopy from any generic immersion of a circle in a two-dimensional plane to an embedding of the circle?  The positive answer to this question was treated as a conjecture (\"{O}stlund conjecture).  In 2014, Hagge and Yazinski disproved the conjecture by showing the first counterexample with a minimal crossing number of 16.  This example is naturally extended to counterexamples with given even minimal crossing numbers more than 14.  
This paper obtains the first counterexample with a minimal crossing number of 15.  This example is naturally extended to counterexamples with given odd minimal crossing numbers more than 13.
%We also obtain infinite families of counterexamples of the \"{O}stlund conjecture, each of which is a triply parameterized curve $P(l, m, n)$.  We study a homotopy described by Reidemeister moves 1 and 3, called (1, 3) homotopy.  We show that $P(l, m, n)$ and $P(l', m', n')$ are (1, 3) homotopic if and only if  $(l, m, n)$ $=$ $(l', m', n')$.
%This study introduces two (1, 3) homotopy invariants for any curve $P$; one is, denoted by $\rii(P)$, the minimum number of Reidemeister moves of type 2 required to obtain the simple closed curve under (1, 3) homotopy and the other is denoted by $c_{\min}(P)$, the minimum number of double points under (1, 3) homotopy.  
%For any nonnegative integer $\mu$, there exits a curve $P$ such that $\rii(P)=\mu$.     
%For a nonnegative integer $\nu$ greater than 14, there exists a curve $P$ such that $c_{\min}(P)$ $=$ $\nu$.   
%Further generalizations of our results are also obtained.     
\end{abstract}  
\section{Introduction}
\begin{figure}[h!]
\includegraphics[width=12cm]{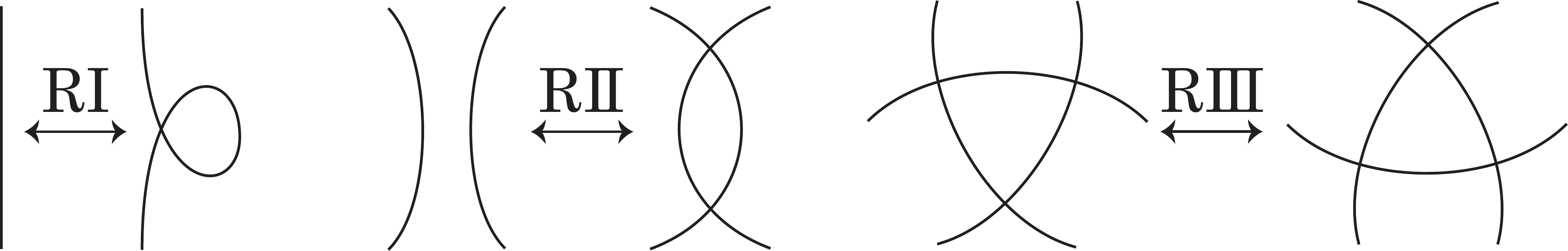}
\caption{$\ri$, $\rii$, and $\riii$}\label{1301}
\end{figure}

A {\it{knot projection}} is the image of a generic immersion of a circle into a $2$-sphere.  Thus, any self-intersection of a knot projection is a transverse double point, which is simply called a {\it{double point}}.  A {\it{trivial knot projection}} is a knot projection with no double points.  
Any two knot projections are related by a finite sequence of Reidemeister moves of types 1, 2, and 3.  Here, Reidemeister moves of types 1, 2, and 3, denoted by $\ri$, $\rii$, and $\riii$, are defined in Fig.~\ref{1301}.
If two knot projections $P$ and $P'$ are related by a finite sequence generated by $\ri$ and $\riii$, then $P$ and $P'$ are {\it{(1, 3) homotopic}}.   The relation becomes an equivalence relation and is called {\it{(1, 3) homotopy}}.  

In 2001 \cite{O}, \"{O}stlund formulated a question as follows: 
\begin{ostlund}
Are Reidemeister moves $\ri$ and $\riii$ sufficient to obtain a homotopy from any generic immersion $S^{1}$ $\to$ $\mathbb{R}^2$ to an embedding?
\end{ostlund}  

In \cite{HY}, Hagge and Yazinski obtained an answer to this question as follows (\cite{HY} treated the positive answer to \"{O}stlund's question as ``his conjecture", and thus, we call it \emph{\"{O}stlund Conjecture}):  
\begin{haggeyazinski}\label{hythm}
A homotopy from the trivial knot projection to $P_{HY}$ that appears as Fig.~\ref{1302} cannot be obtained by a finite sequence generated by Reidemeister moves $\ri$ and $\riii$.     
\begin{figure}[h!]
\includegraphics[width=2cm]{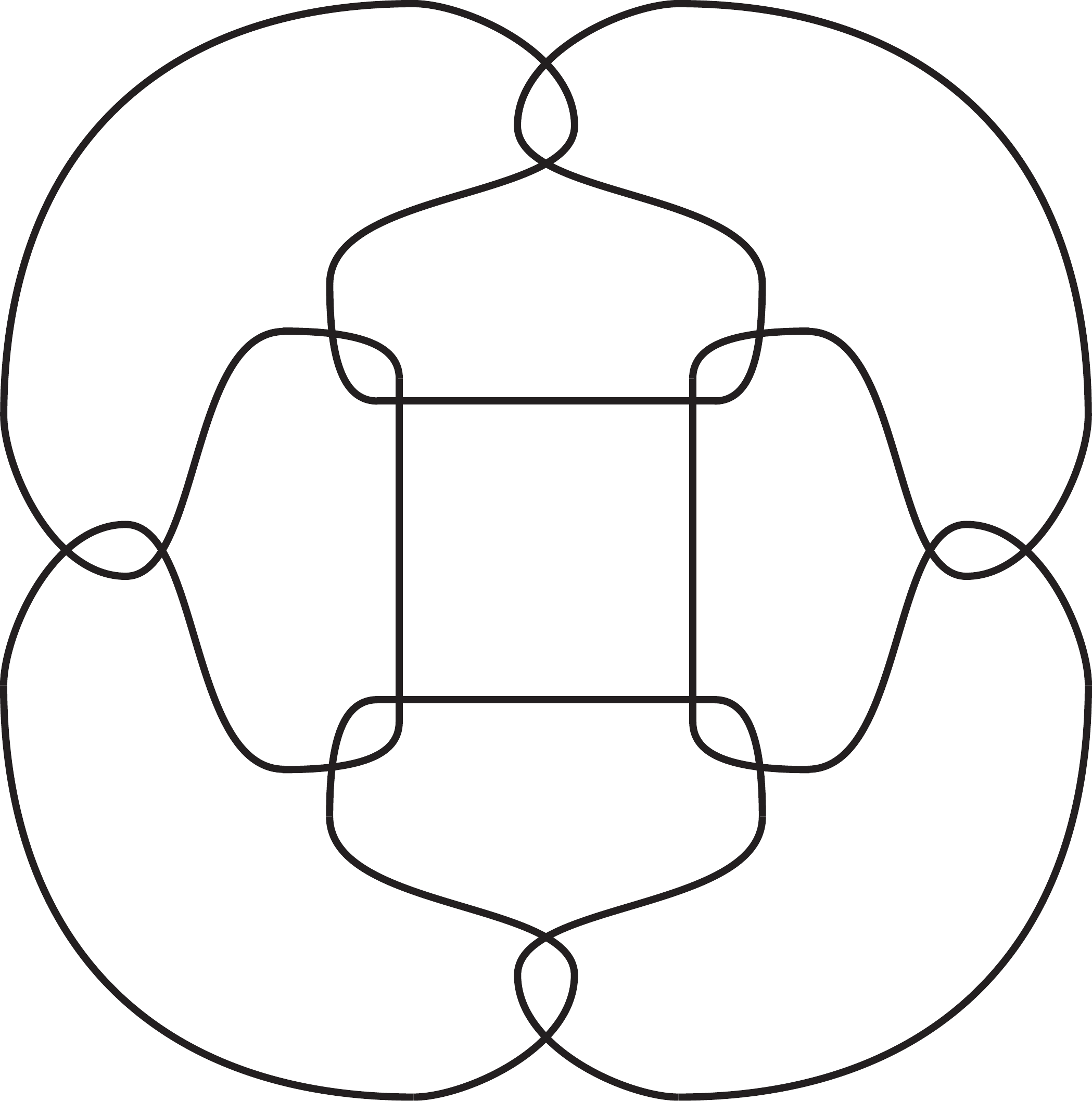}
\caption{Hagge-Yazinski's example $P_{HY}$}\label{1302}
\end{figure}
\end{haggeyazinski}    

For a given knot projection $P$, let $c_{\min}(P)$ $=$ $\min\{$ the number of double points of $P'~|~P'$ and $P$ are related by a finite sequence generated by $\ri$ and $\riii$ $\}$.  We call $c_{\min}(P)$ the \emph{minimal crossing number} of $P$.  

As a step further, it is easy to find a vast generalization of the Hagge-Yazinski Theorem; there exists an infinite family of counterexamples of the \"{O}stlund Conjecture, i.e., non-trivial knot projections under (1, 3) homotopy (Remark~\ref{remark_hy}).  
However, every knot projection of this family has an even minimal crossing number.  

In this paper, we find a counterexample of the conjecture with an odd minimal crossing number.  
We show the non-triviality of the example by using Hagge-Yazinski techniques \cite{HY}.  
This example is naturally extended to an infinite family of knot projections, each of which is a knot projection with a given odd minimal crossing number.   
%We also introduce two invariants, $\rii(P)$ and $c_{\min}(P)$, under (1, 3) homotopy.  We show that for any nonnegative integer $\mu$, there exists a knot projection $P$ such that $\rii(P)=\mu$ (Theorem~\ref{main1}~(2)); here, by definition, for the trivial knot projection $U$, $\rii(U)=0$.  We also show that for any integer $\nu$ ($\ge 15$), there exists a knot projection $P$ such that $c_{\min}(P)=\nu$ (Theorem~\ref{main1}~(3) and Theorem~\ref{main2}~(\ref{main2_3})).  
%\begin{theorem}\label{main1}
Similar to Remark~\ref{remark_hy}, Theorem~\ref{main2} naturally implies a knot projection $P$ such that $c_{\min}(P)$ $=$ $2i+13$ ($i \in \mathbb{Z}_{> 0}$).
\begin{theorem}\label{main2}
The knot projection, shown in Fig.~\ref{1320}, and the trivial knot projection cannot be (1, 3) homotopic.  
\begin{figure}[h!]
\includegraphics[width=2.5cm]{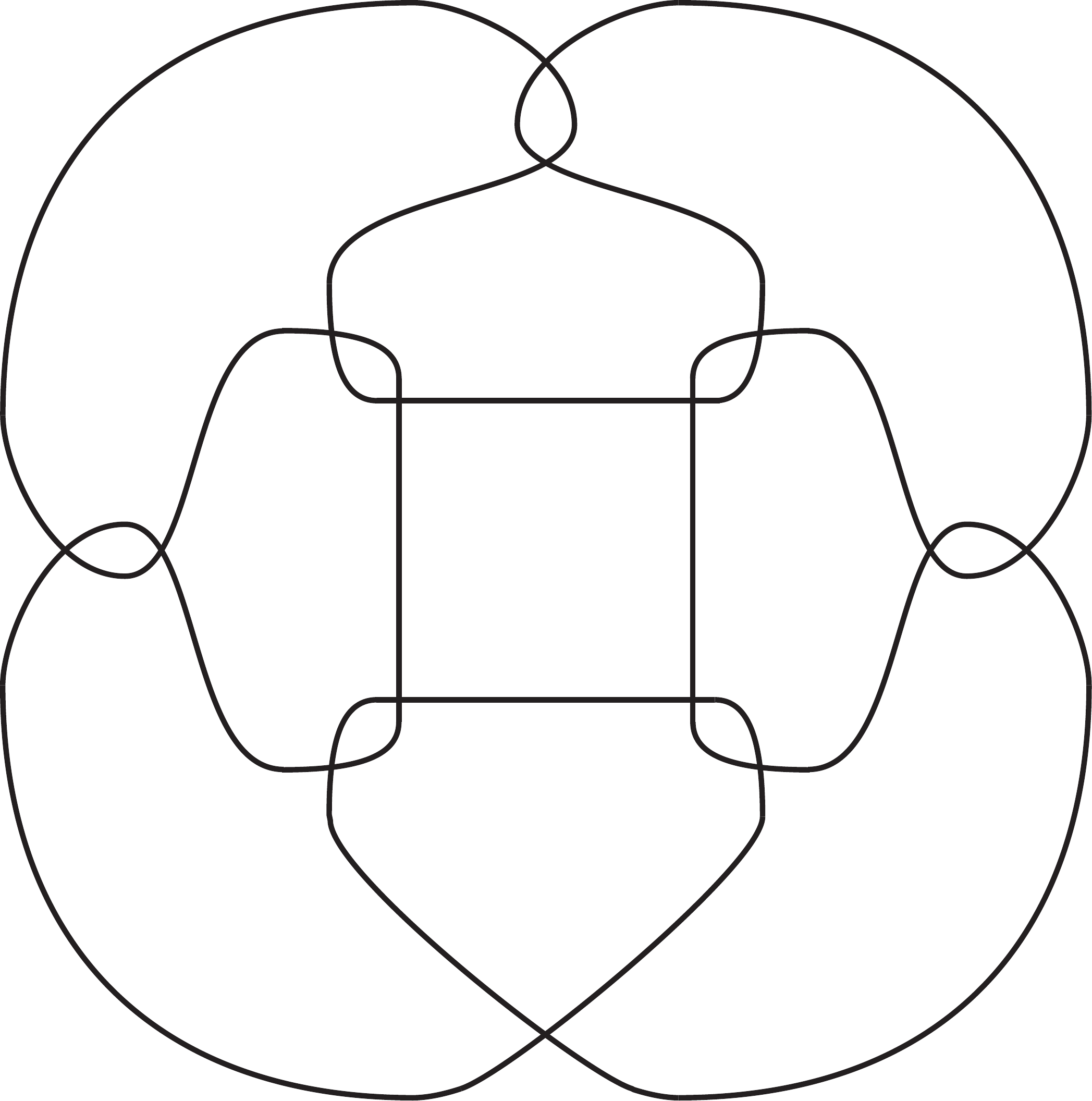}
\caption{A knot projection}\label{1320}
\end{figure}
\end{theorem}

\section{Proof of Theorem~\ref{main2}}\label{sec_proof}    
\begin{figure}[h!]
\includegraphics[width=10cm]{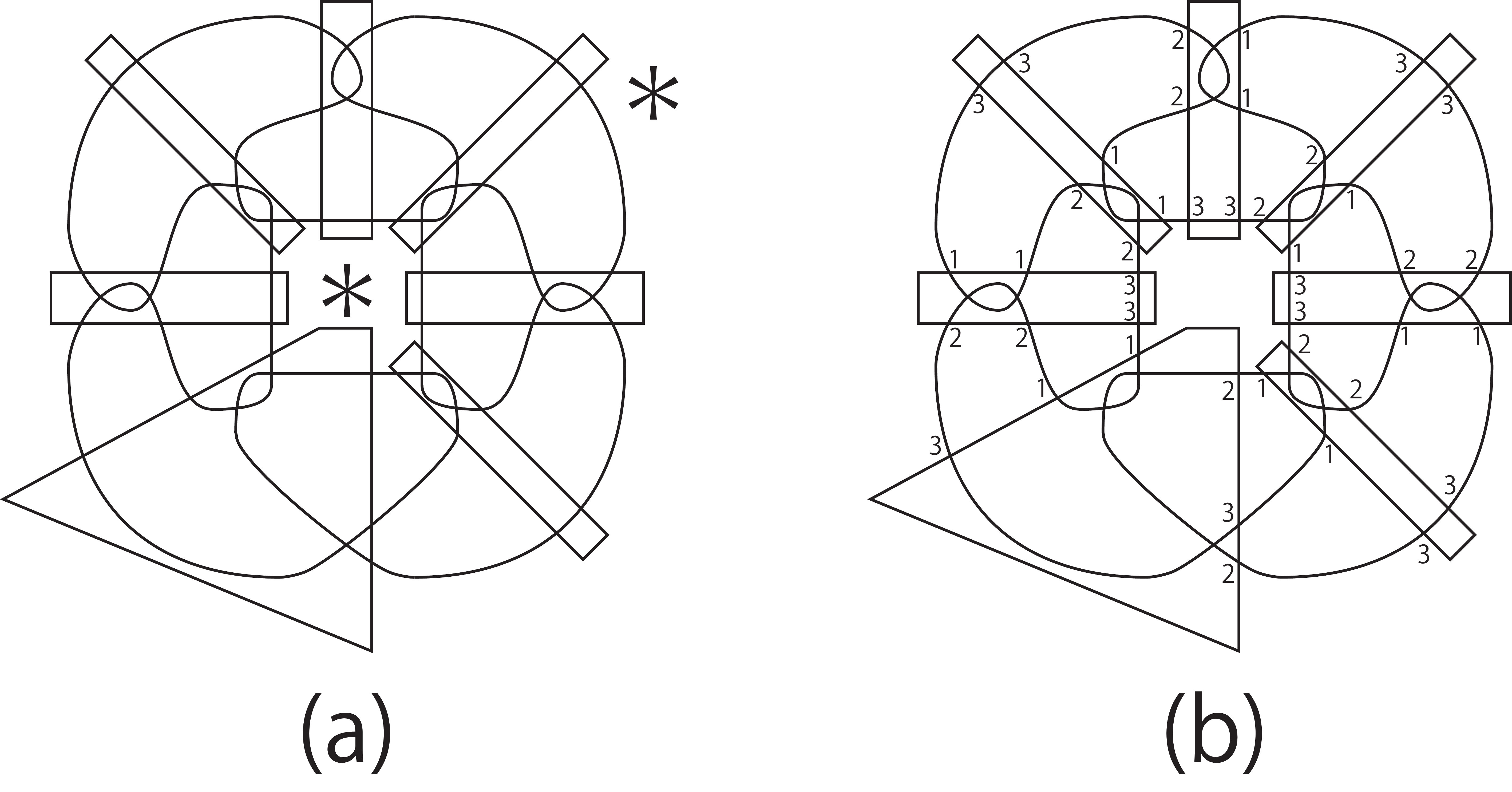}
\caption{(a): seven boxes and two starred polygons and (b) the configuration of strands}\label{1317}
\end{figure}
As shown in Fig.~\ref{1317}, any $P_0$ can be decomposed into seven boxes and arcs with no double points.  Each box is equivalent to $[0, 1] \times [0, 1]$ under sphere isotopy and in each box, there is a part of a knot projection, called (3, 3)-{\it{tangle}} (cf.~Definition~\ref{def_tangle}).  
\begin{definition}[($s, t$)-tangle]\label{def_tangle}
Let $s$ and $t$ be positive integers such that $s+t$ is even.  An unoriented ($s, t$)-tangle is the image of a generic immersion of $(s+t)/2$ arcs into $[0, 1] \times [0, 1]$ such that:
\begin{itemize}
\item The boundary points of the arcs map bijectively to $s+t$ points \[ \{1, 2, \dots s \} \times \{ 1 \}, \{ 1, 2, \dots, t \} \times \{ 0 \}.\]  
\item Near the endpoints, the arcs are perpendicular to the boundary $[0, 1]$.  
\end{itemize}
An image of the map of a single arc is called a {\it{strand}}.  
\end{definition}
\begin{definition}[polygon, $r$-gon]
Let $P$ be a nontrivial knot projection and let $r$ be a positive integer.  For each connected component $D$ of $S^2 \setminus P$, $\partial D$ consists of a circle with double point(s).  Then, $D \cup \partial D$ is called a \emph{polygon}.  If $\partial D$ has exactly $r$ double point(s), $D \cup \partial D$ is also called $r$-\emph{gon}. 
\end{definition} 

Each of the seven boxes satisfies the following four conditions, called the {\it{box property}} consisting of Rules (1)--(3):
\begin{enumerate}
\item There exist exactly two polygons, called {\it{starred polygons}}, each having at least four sides partially outside boxes and each completely containing seven boxes' sides with no intersections with the knot projection.  \label{r3}
\item Choose and fix one of the starred polygons.  Because we fix our gazing direction based on the selected starred polygon, we can define the {\it{left-side}} and {\it{right-side}} of each box.  Strand~1 (resp. $2$) is a strand that begins and ends on the left-side (resp.~right-side) in each box.  Strand~3 has one endpoint on the left-side and another endpoint on the right-side of each box.  
Strands~1 and $2$; for each pair, strands~1 and 2 intersect at exactly $2$ double points.     
\label{r1}

\item No double points are placed outside the boxes.  Each endpoint of strand~3 on the left-side is connected to the endpoint of strand~2 on the right-side of the adjacent box.  Each endpoint of a strand~3 on the right-side is connected to the endpoint of  strand~1 on the left-side of the adjacent box.  The left pairs of strands~1 and 2 are connected to each other.   \label{r2} 
%\item There exists at least one box in which strand~3 intersects either strand~1 or strand~2.    A box satisfying this condition is called a {\it{special box}}.   \label{r4}
\end{enumerate}

%\noindent $\bullet$ {\it{Proof of $\rii(P(l, m, n)) \ge 1$.}}
We consider an inductive proof with respect to the number of Reidemeister moves, $\ri$ or $\riii$, applied to $P_0$.  This induction proves the existence of seven boxes satisfying the box property, which implies that $P_0$ cannot be (1, 3) homotopic to the trivial knot projection (cf.~Rule~(\ref{r1})).  Now we prove the claim.  Let $n$ be the number of Reidemeister moves applied to $P_0$.  

First, for the case where $n=0$, there exist seven boxes such that the box property holds.  Second, we assume the existence of seven boxes such that the box property holds for $n$ and prove the $(n+1)$th case.  Let 1a (resp.~1b) be $\ri$ increasing (resp.~decreasing) the number of double points.  
In the following, we show the $(n+1)$th claim for each Reidemeister move within a box and outside or partially outside the boxes.  From the assumption (Rule~(\ref{r3})), the definition of two starred polygons for the $n$th case are well-defined.  
%This is because from Rule~(\ref{r1}) and Rule~(\ref{r2}) for the $n$th case we can determine the connections of all boxes; thus it can be seen that the corresponding polygon still has at least four sides partially outside boxes after the applications of $n$ Reidemeister moves to the initial starred polygon.  
We call each starred polygon an $n$th {\it{starred polygon}}.  
\subsection{$\ri$ or $\riii$ occurring within a box}
Let us consider a single move (RI or RI\!I\!I) occurring entirely within one of the boxes.  This move fixes the endpoints of the strands, and thus Rule~(\ref{r1}) and Rule~(\ref{r2}) are satisfied.  
That is, if we assume that the $n$th case is true, we can define the same number of boxes satisfying Rules~(\ref{r1}) and (\ref{r2}) for the ($n+1$)th case after the application of $\ri$ or $\riii$ within a box.     

\begin{figure}[h!]
\includegraphics[width=7cm]{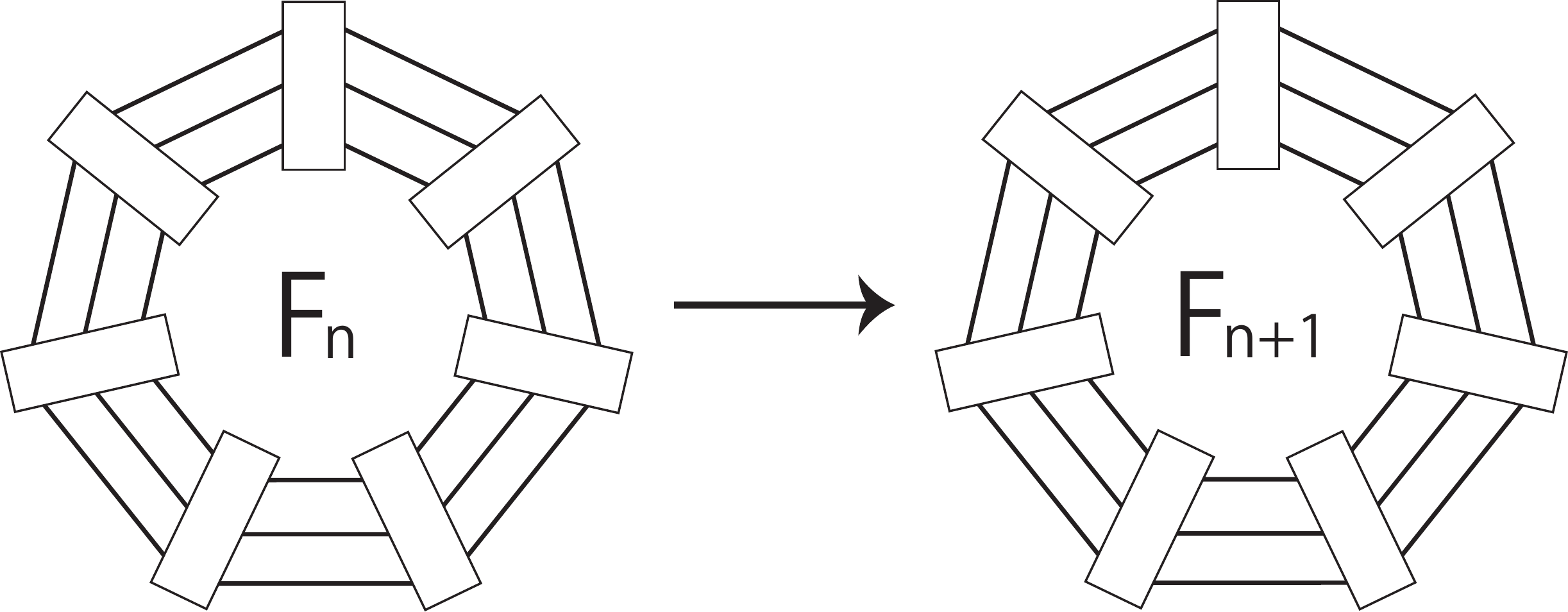}
\caption{$F_n$ and $F_{n+1}$}\label{1306_5}
\end{figure}

From the above condition, we can prove Rule~(\ref{r3}) using Rules~(\ref{r1}) and (\ref{r2}) even if the order of the endpoints of the strands is changed for a side.  
Let $P_n$ and $P_{n+1}$ be knot projections.  
Let $F_{n+1}$ be a polygon of $P_{n+1}$ corresponding to an $n$th starred polygon $F_n$ of $P_n$ where $P_{n+1}$ is obtained by applying a single Reidemeister move to $P_{n}$ within a box (Fig.~\ref{1306_5}).  

Suppose that $F_{n+1}$ is an $r$-gon.  Polygon $F_{n+1}$ faces no successive strands of type~3 with no double points, within two successive boxes.  This is because strand~3 is connected to the endpoints of strands~1 or 2 and there is no double point outside the boxes as per Rule~(\ref{r2}).  
Here, note that from Rule~(\ref{r1}), any strand with endpoints on left- and right-sides is strand~3.  
Thus, proceeding along the boundary of $F_{n+1}$, we encounter boxes containing at least four vertices of $F_{n+1}$ that are inside the boxes, which implies $r \ge 4$.    
Thus, Rule~(\ref{r3}) also holds for the $(n+1)$th case.  

In conclusion, if the $(n+1)$th move is within a box, there exist seven boxes such that the box property holds.  
\subsection{1a not occurring within a box}
\begin{figure}
\includegraphics[width=10cm]{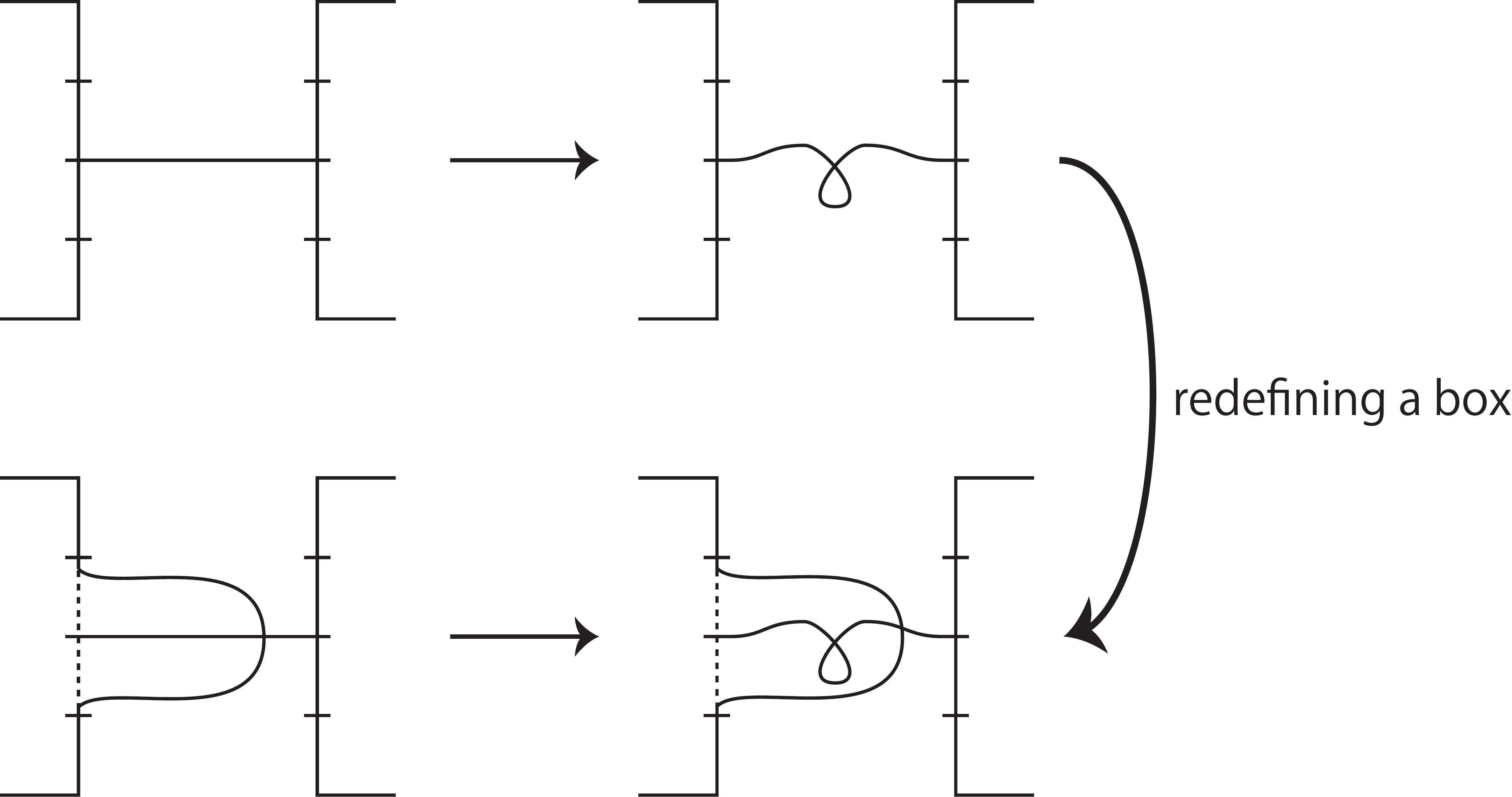}
\caption{Redefining a box}\label{1307}
\end{figure}

Note that we assume the existence of boxes such that the box property holds for the $n$th case.  
If the $(n+1)$th move is 1a and is outside the boxes, we can redefine the box by sphere isotopy, as shown in Fig.~\ref{1307}.  If this 1a is not completely outside the boxes, we can use a similar modification.   
\subsection{1b not occurring within a box}
If the $(n+1)$th move is 1b and not occurring within a box, there exists a $1$-gon to be removed for the $n$th case.  The two possibilities of the $1$-gon are considered as follows.
\begin{itemize}
\item From Fig.~\ref{1308},   
\begin{figure}
\includegraphics[width=3cm]{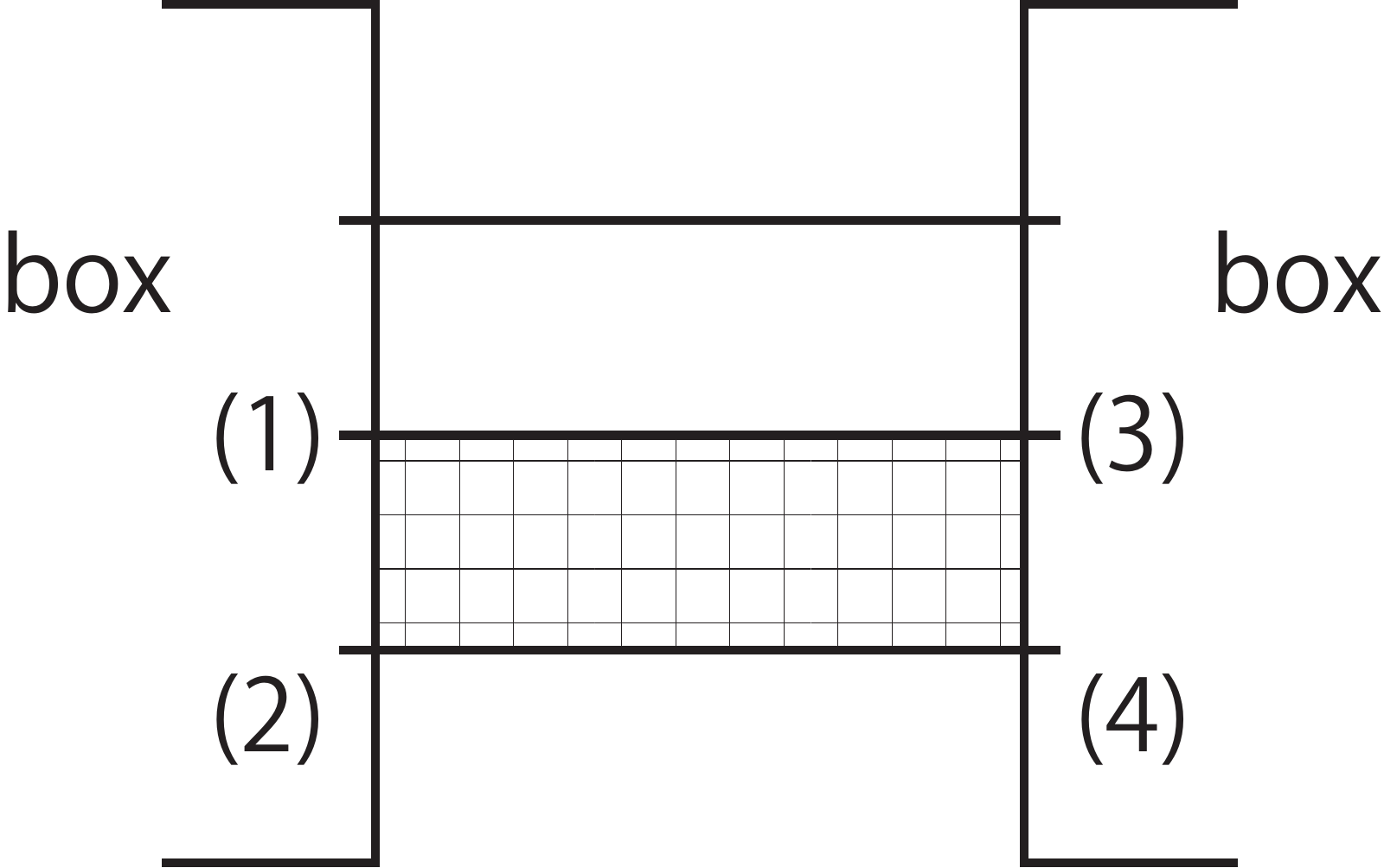}
\caption{Each of endpoints (1), (2), (3), and (4) connects to strand~1, 2, or 3.}\label{1308}
\end{figure}
we can see the $1$-gon of 1b containing a region having two sides outside the boxes.  If one of the two sides is directly connected to either strand~1 or 2, then the $1$-gon has at least two double points, which is a contradiction (there is no $1$-gon with two double points).  If one of the two sides is directly connected to strand~3 that is connected strand~1 or 2 in the adjacent box, then this also implies a contradiction similar to the case above. 
\item For $1$-gon of 1b containing an $n$th starred polygon, by the induction assumption, it has at least $r$ ($\ge 4$) vertices, which is a contradiction.  
\end{itemize}  
\subsection{$\riii$ not occurring within a box}
\begin{figure}[h!]
\includegraphics[width=5cm]{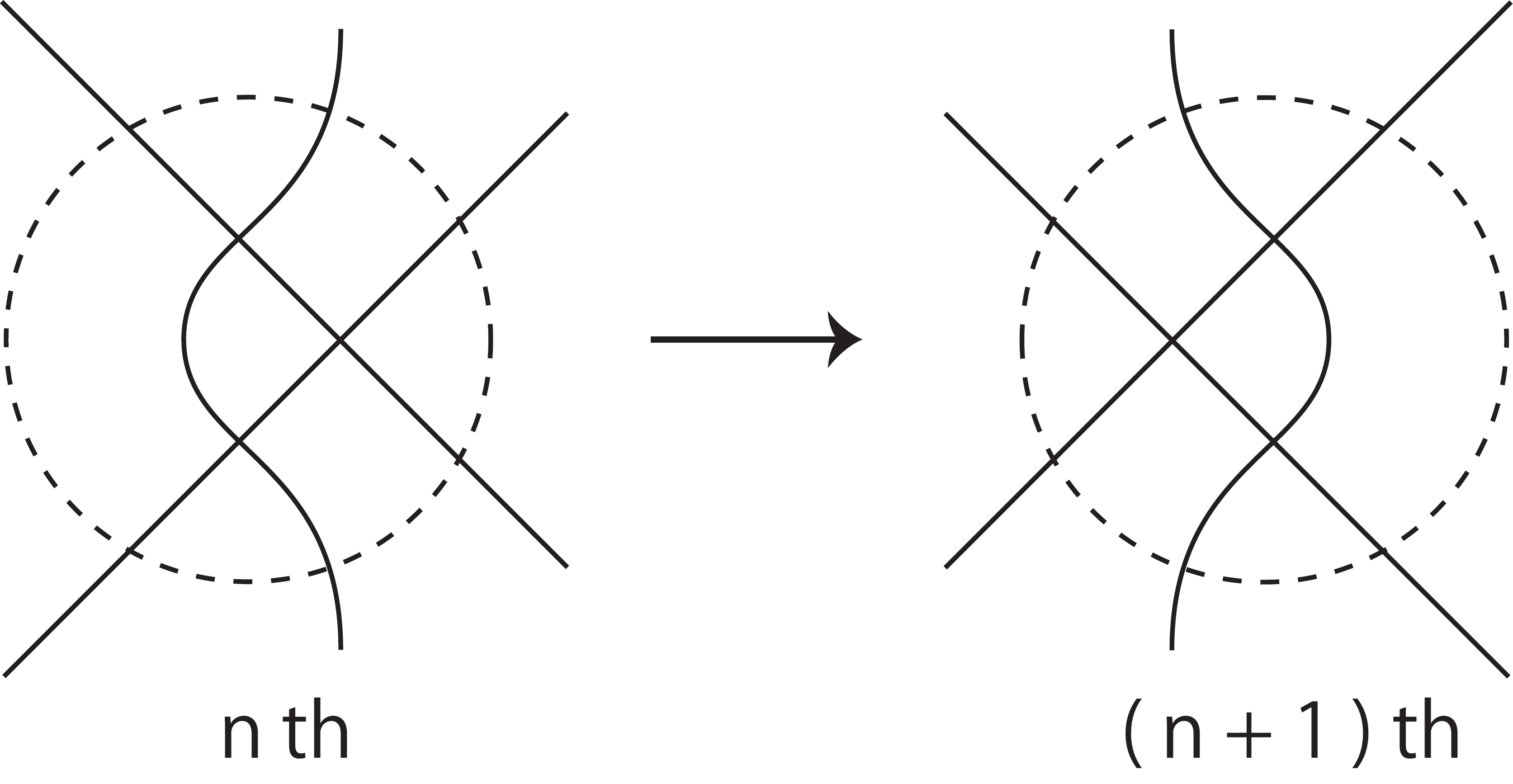}
\caption{$3$-gon $T_n$ for the $n$th case with respect to a single~$\riii$}\label{1309}
\end{figure}

\begin{figure}[h!]
\includegraphics[width=3cm]{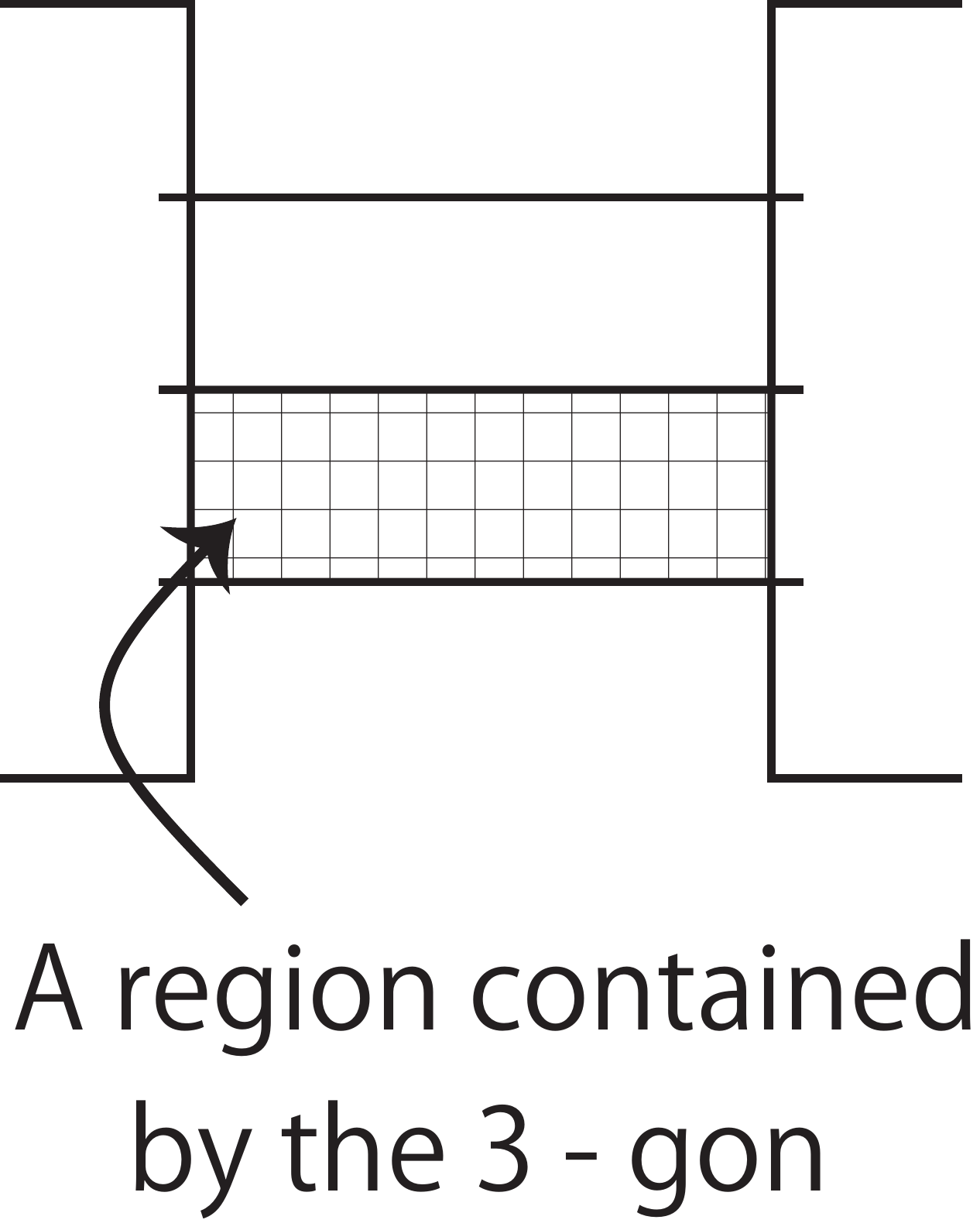}
\caption{A region contained by the $3$-gon appearing in the $n$th case}\label{1310}
\end{figure}

If the $(n+1)$th move is a single $\riii$, we focus on the $3$-gon $T_n$ for the $n$th case with respect to the single $\riii$, as shown in Fig.~\ref{1309}.  
If $T_n$ is not in a box, it contains a region having two sides outside the box, as shown in Fig.~\ref{1310}.  
This is because $T_n$ cannot contain a starred polygon ($=$ $r$-gon, $r$ $\ge 4$) by the induction assumption.

\begin{figure}[h!]
\includegraphics[width=3cm]{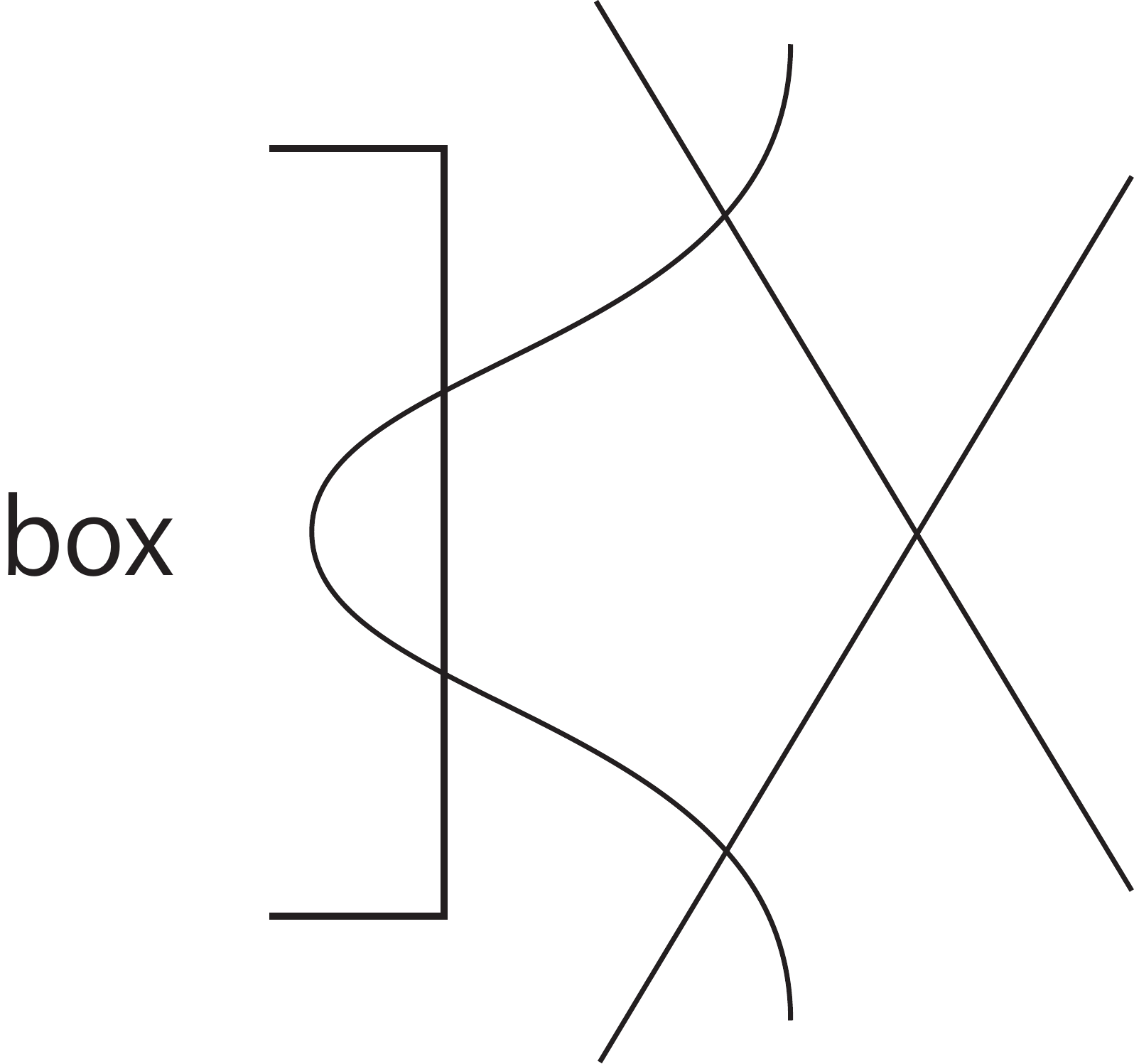}
\caption{Impossible case}\label{1311}
\end{figure}

Further, if $T_n$ is not inside a box, at least one vertex of $T_n$ is inside a box.  If not, there is an arc inside a box, as shown in Fig.~\ref{1311}, which starts and ends on the same side.   
However, by the induction assumption and Rule~(\ref{r1}) of the box property for the $n$th case, the arc should be either strand~1 or 2, which is a contradiction because an arc has no intersection whereas both strand~1 and 2 have intersections.  
  
\begin{figure}[h!]
\includegraphics[width=5cm]{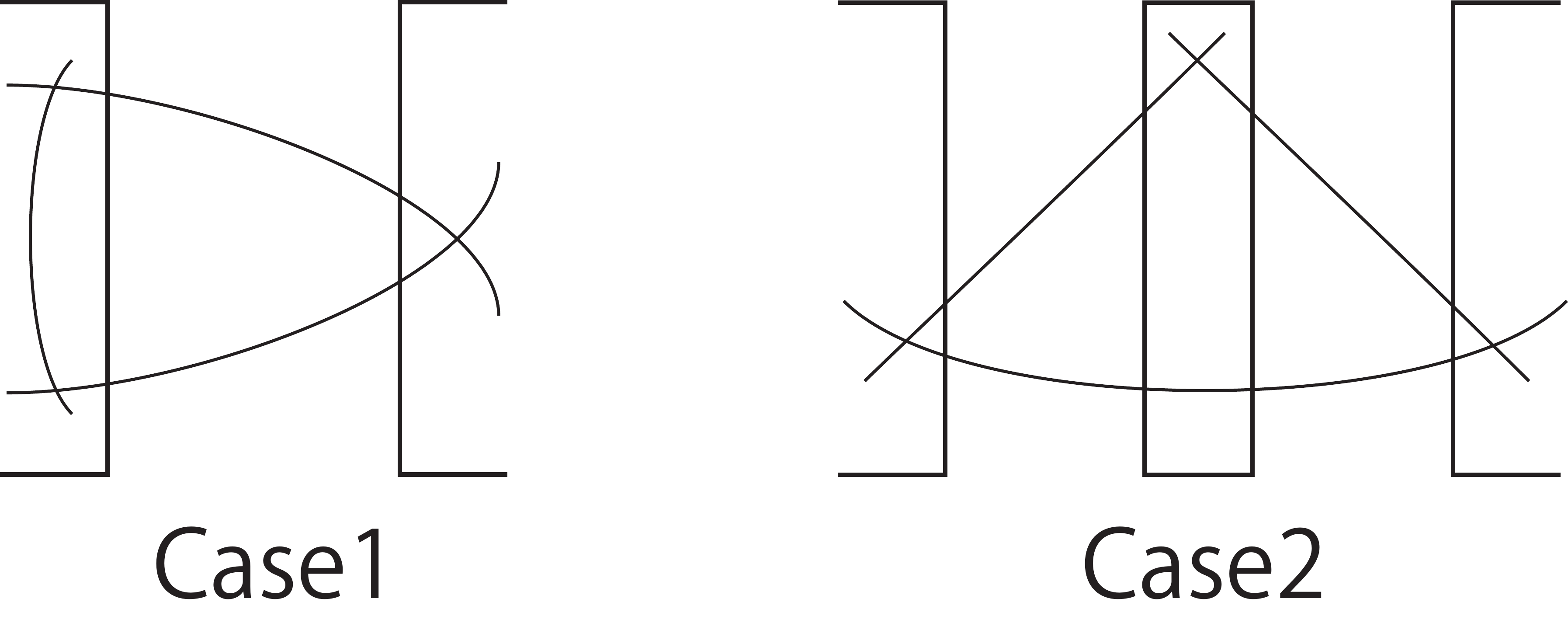}
\caption{Case~1 (left) and Case~2 (right)}\label{1312}
\end{figure}

\begin{figure}[h!]
\includegraphics[width=8cm]{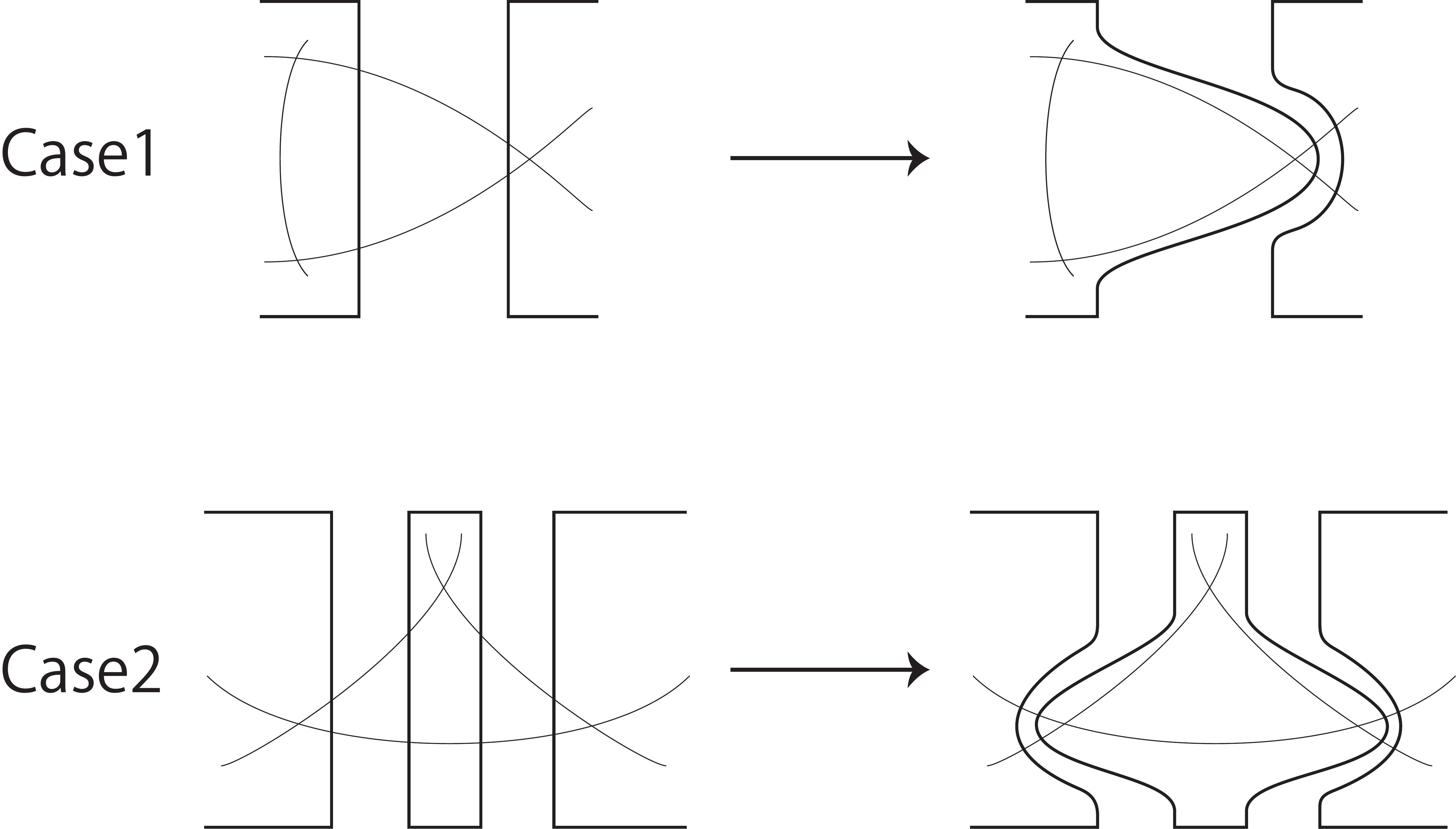}
\caption{Redefining box(es) for Case 1 (upper) and Case 2 (lower)}\label{1313}
\end{figure}

By the induction assumption, there is no double point outside the boxes for the $n$th case.  Thus, if $T_n$ is not inside a box, there are exactly two cases of local arrangements of boxes.  Case~1 is when one vertex of $T_n$ is inside a box and two vertices of $T_n$ are in another box.  Case~2 is when all three vertices are in different boxes (Fig.~\ref{1312}).    
As a result, for Case~1 (resp.~Case~2), by sphere isotopy, we can redefine one box (resp.~two boxes), as shown in Fig.~\ref{1313}.  

Here, it is easy to see that this redefinition satisfies Rules~(\ref{r3})--(\ref{r2}).  
%In the following, we show the existence of special box (Rule~(\ref{r4})).    For Case~1, if a box containing two vertices of a $3$-gon with respect to RI\!I\!I is not a special box, retake the box, as shown in Fig.~\ref{1313}.  In this case, there is a possibility that the special box is moved to the next box.  However, there exists at least one special box, which implies Rule~(\ref{r4}).    If a box containing two vertices is the special box, redefine the boxes, as shown in Fig.~\ref{1313} (Case~1).  In this case, this special box for the $n$th case is a box satisfying the box property.  In particular, the box satisfies the box property and can be called a special box for $(n+1)$th case.  There exists at least one special box, which implies Rule~(\ref{r4}).    For Case~2, redefining of the boxes as shown in Fig.~\ref{1313} (Case~1).  In this case, there is a possibility that the special box is moved to the next box (maybe twice), but there exists at least one special box, which implies Rule~(\ref{r4}).     
Thus, the $(n+1)$th case holds and this completes the proof.  
$\hfill\Box$

\begin{remark}\label{remark_hy}
For positive integers $l \ge 1$, $m \ge 1$, and $n \ge 4$, let $P(l, m, n)$ be a knot projection, as defined in Fig.~\ref{1303}.  It is easy to see that $P_{HY}$ is extended to $P(l, m, n)$.  
\begin{figure}[h!]
\includegraphics[width=12cm]{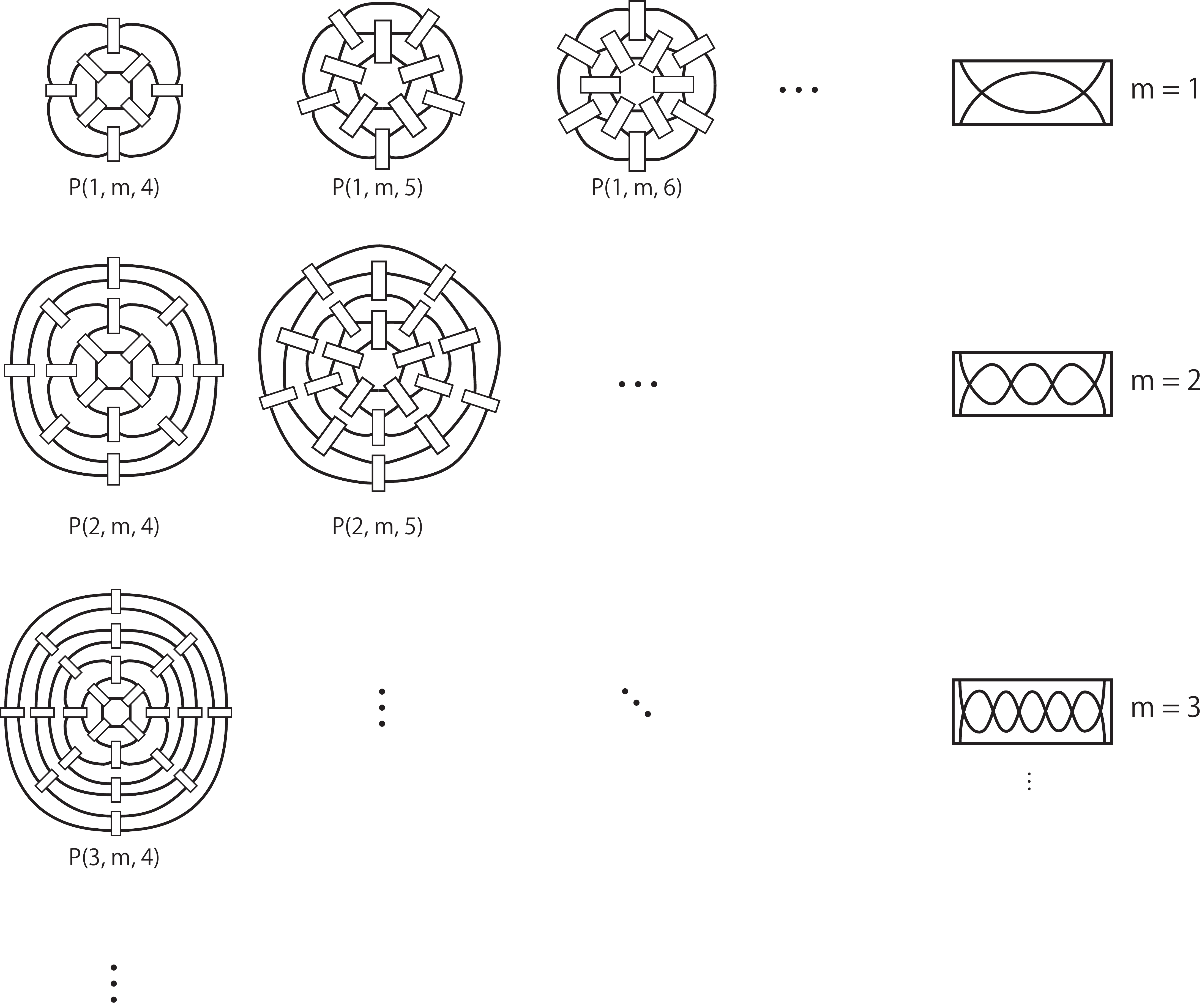}
\caption{$P(l, m, n)$}\label{1303}
\end{figure}
It is also elementary to find that for any positive integer $i$, there exists a knot projection with $c_{\min}(P)$ $=$ $2i+14$ (see Fig.~\ref{1304}).   
%In particular, if $i \neq i'$, $P_i$ and $P_{i'}$ are not (1, 3) homotopic.   
\begin{figure}[h!]
\includegraphics[width=7cm]{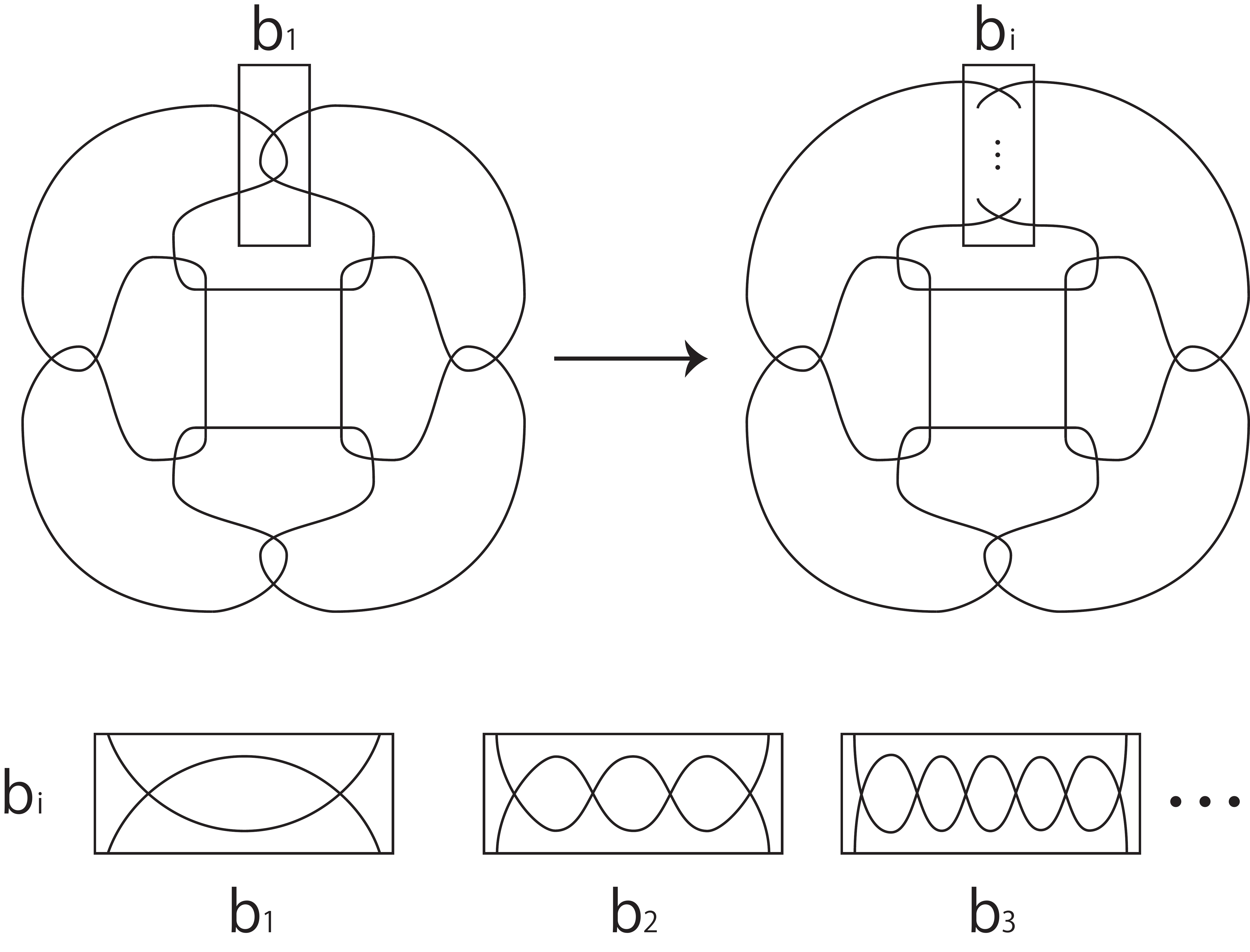}
\caption{A knot projection $P$ with $c_{\min}(P)=2i+14$}\label{1304}
\end{figure}
\end{remark}
%\begin{remark}
%A homotopy from every knot projection with at most 7 double points to the trivial knot projection is obtained by a finite sequence generated by RI and RI\!I\!I.   
%\end{remark}

\section*{Acknowledgements} 
The authors would like to thank Professor Takahiro Kitayama for his useful comments at the conference: Intelligence of Low-dimensional Topology, RIMS, Kyoto University, held in May 2015.  
We also thank Professor Tomoo Yokoyama for his comments on our earlier work at the seminar: KUAMS, ACCA-JP, Sakajo CREST seminar, Kyoto University, held in November 2014.  
The work of N.~Ito was partially supported by the JSPS Japanese-German Graduate Externship and a Waseda University Grant for Special Research Projects (Project number: 2015K-342).  N.~Ito is a project researcher of Grant-in-Aid for Scientific Research (S) 24224002 (April 2016--).

\end{document}